\documentclass[10pt]{article}

\usepackage{amsfonts,amssymb,amsmath,fontenc}
\usepackage{latexsym,wasysym,mathrsfs}
\usepackage{graphics,enumerate}
\usepackage{color,graphicx}
\usepackage[pdftex]{hyperref}
\usepackage[left=3cm,top=3cm,right=3cm,bottom=3cm]{geometry}

\definecolor{darkred}{rgb}{.0,.0,.80}
\definecolor{darkgreen}{rgb}{.0,0,.80}
\definecolor{darkredd}{rgb}{0,0,.80}
\hypersetup{
 colorlinks,%
       citecolor=darkgreen,%
    filecolor=magenta,%
    linkcolor=darkredd,%
    urlcolor=darkred,%
}

\DeclareMathAlphabet{\mathpzc}{OT1}{pzc}{m}{it}

\title{The Cauchy problem for the DMKP equation}

\date{}

\author{\textbf{Amin Esfahani}\footnote{Address: School of Mathematics and Computer Science, Damghan University, Damghan, Postal Code 36716-41167, Iran.\, E-mail: amin@impa.br, esfahani@du.ac.ir}}

\date{}

\begin{document}
\maketitle

\begin{abstract}
In this work, we study the dissipation-modified Kadomtsev-Petviashvili equation in two  space-dimensional case.
We establish  that the Cauchy problem for this equation is locally well-posed in anisotropic Sobolev spaces. We show in some sense that our result is sharp. We also prove the global well-posedness for this equation under  suitable conditions.
\end{abstract}

\numberwithin{equation}{section}

\newtheorem{theorem}{{ THEOREM}}[section]
\newtheorem{lemma}[theorem]{{\sc \textbf{LEMMA}}}
\newtheorem{defn}[theorem]{{\sc \textbf{DEFINITION}}}
\newtheorem{proposition}{\begin{small}PROPOSITION\end{small}}
\newtheorem{remark}[theorem]{{\sc \textbf{REMARK}}}
\newtheorem{corollary}[theorem]{{\sc \textbf{Corollary}}}
\newcommand{\rr}{{\mathbb R}}
\newcommand{\dl}{\Delta}
\newcommand{\T}{{\mathbb T}}
\newcommand{\C}{{\mathbb C}}
\newcommand{\N}{{\mathbb N}}
\newcommand{\lam}{{\lambda}}
\newcommand{\ff}{\varphi}
\newcommand{\ti}{\widetilde}
\newcommand{\what}{\widehat}
\newcommand{\Z}{{\mathbb Z}}
\newcommand{\om}{\omega}
\newcommand{\pf}{\noindent\textbf{Proof.}\quad}
\newcommand{\fim}{\hfill$\square$\\ \\}
\newcommand{\dd}{{\rm d}}
\newcommand{\ee}{{\rm e}}
\newcommand{\ii}{{\rm i}}
\newcommand{\es}{^{\APLstar}}
\newcommand{\sech}{\textmd{sech}}
\newcommand{\vep}{\varepsilon}
\newcommand{\vr}{\varrho}
\newcommand{\vrx}{\varrho(\xi)}

\newcommand{\cchi}{\mathscr{R}(\zeta,\zeta_1)}
\newcommand{\M}{\mathscr{M}(\xi,\xi_1)}
\newcommand{\K}{\mathscr{K}(\zeta,\zeta_1,t)}

\section{Introduction}
In this work, we consider the initial value problem for  the dissipation-modified Kadomtsev-Petviashvili (DMKP) equation
\begin{equation}\label{dmkp}
\left\{\begin{array}{ll}
{\left(u_t+u_{xxx}+uu_x+\alpha(u_{xx}+u_{xxxx})+\beta(u^2)_{xx}\right)}_x+\vep u_{yy}=0,\qquad (x,y)\in\rr^2,\,\,t\geq0,\\\\
u(x,y,0)=\ff(x,y),
\end{array}\right.\end{equation}
where $\alpha>0$ and $\beta$ are real constants and $\vep=\pm1$. The DMKP equation \eqref{dmkp} arises in studying spontaneous generation of long waves in the presence of a conservation law in isotropic systems (e.g., B\'{e}nard-Marangoni waves), near the instability threshold \cite{bar,Gvelarde,nvelarde}. In \cite{prep,prep1}, the author has also investigated another version of \eqref{dmkp}.

Equation \eqref{dmkp} is also a natural two-dimensional version of the KdV-Kuramoto-Sivashinsky (KdV-KS) equation
\begin{equation}\label{kdvks}
u_t+u_{xxx}+uu_x+\alpha(u_{xx}+u_{xxxx})=0,
\end{equation}
which arises in interesting physical situations, for example as
a model for long waves on a viscous fluid flowing down an inclined plane \cite{topperkawahara} and
to derive drift waves in a plasma \cite{cohen}.

The DMKP equation, when $\beta=0$, is a dissipative version of the Kadomtsev-Petviashvili (KP) equation
\begin{equation}\label{kp}
(u_t+u_{xxx}+uu_x)_x+\vep u_{yy}=0,
\end{equation}
which is universal model for nearly one directional weakly nonlinear dispersive waves with weak transverse effects. The KP equation, in turn,  is a  two-dimensional extension of the KdV equation
\begin{equation}\label{kdv}
u_t+u_{xxx}+uu_x=0.
\end{equation}
Our principal aim here is to study the local well-posedness for the initial value problem associated to the DMKP equation in the anisotropic Sobolev spaces $H^{s_1,s_2}(\rr^2)$, $s_1,s_2\in\rr$.

In the past decades, Bourgain developed a new method, clarified by Ginibre in \cite{ginibre}, for the study of the Cauchy problem for nonlinear dispersive equations. This method was further successfully applied to Schr\"{o}dinger, KdV  and KPII equations (cf. \cite{bourgain,bourgain-2,bourgain-3,kpv1,kpv2}).  The original Bourgain method makes extensive use of the
Strichartz inequalities in order to derive the bilinear estimates corresponding to the nonlinearity. On
the other hand, Kenig et al. \cite{kpv1,kpv2} simplified Bourgain's proof and improved the bilinear estimates using
only elementary techniques, such as Cauchy-Schwartz inequality and simple calculus inequalities (see also \cite{GTV,tao}).

It was also shown by Molinet and Ribaud \cite{molinetribaud-1,molinetribaud,molinetribaud-2} that the Bourgain spaces can be used to study the Cauchy problems associated to semi-linear equations with a linear part containing both dispersive and dissipative terms; and consequently this applies to the KdV-Burgers (KdVB) equation
\begin{equation}
u_t+u_{xxx}+uu_x=u_{xx}
\end{equation}
and the  Kadomtsev-Petviashvili-Burgers (KPB) equation
\begin{equation}\label{kpb}
(u_t+u_{xxx}+uu_x-u_{xx})_x+\vep u_{yy}=0.
\end{equation}
By introducing a Bourgain space associated to the usual KP equation, related only to the
dispersive part of the linear symbol of \eqref{kpb}, Molinet and Ribaud \cite{molinetribaud}  proved global existence
for the Cauchy problem associated to the KPB equation \eqref{kpb}, by using Strichartz-type estimates for the KP equation injected into the framework of Bourgain spaces. More precisely, authors in \cite{molinetribaud} showed the KPB-I equation ($\vep=-1$) is locally well-posed in $H^{s_1,s_2}(\rr^2)$ if $s_1>0$ and $s_2\geq0$;  and the KPB-II equation ($\vep=1$) is locally well-posed in $H^s(\rr^2)$ if $s\geq0$. The global well-posedness followed by means of a priori estimates.

Recently Kojok in \cite{kojok} obtained a sharp result by proving that the  KPB-II equation is globally well-posed in
$H^{s_1,s_2}(\rr^2)$ for $s_1>-1/2$ and $s_2\geq0$.

In this paper, we will apply the ideas of  \cite{molinetribaud-1,molinetribaud,molinetribaud-2} and introduce a Bourgain-type space associated to
the KP equation. This space is in fact the intersection of the space introduced in \cite{bourgain-3} and of a Sobolev space. The advantage of this space is that it contains both the dissipative and dispersive parts of the linear symbol of \eqref{dmkp}. Next we  establish the local existence for \eqref{dmkp} with initial value $\ff\in H^{s_1,s_2} (\rr^2)$ when $s_1 > -1/2$ and $s_2\geq0$; and we also show that the Cauchy problem \eqref{dmkp} is globally well-posed in $\ff\in H^{s_1,s_2} (\rr^2)$ if $\beta=0$, $s_1 > -1$ and $s_2\geq0$. We prove also that our local existence theorem is optimal by constructing a counterexample showing that the application $\ff\to u$ from $H^{s_1,s_2}(\rr^2)$ to $C([0,T ];H^{s_1,s_2}(\rr^2))$ cannot be regular
for $s_1 < -1/2$ and $s_2 = 0$.

This existence result, in some sense, is quite surprising. There is no difference in the existence result
for $\vep=\pm1$ in \eqref{dmkp}. However,
despising the dissipation terms in \eqref{dmkp}, we obtain the KP equation \eqref{kp}, where the KP-I and KP-II models are quite distinct.

This paper is organized as follows. In Section \ref{2}, we introduce some notations and our main results. In Section \ref{3}, we derive linear estimates and some smoothing properties for the operator arising from \eqref{dmkp}
in the Bourgain spaces (Lemma \ref{lem2}). Section \ref{4}  is devoted to establish bilinear estimates by using Strichartz-type estimates for the KP equation.
In Section \ref{5}, using bilinear estimates,
a standard fixed point argument and some smoothing properties, we prove uniqueness and
local existence of the solution of \eqref{dmkp} in anisotropic Sobolev space $H^{s_1,s_2}(\rr^2)$ with $s_1 > -1/2$ and $s_2\geq 0$; and global existence of the solution of \eqref{dmkp}, with $\beta=0$, in  $H^{s_1,s_2}(\rr^2)$ with $s_1 > -1$
and $s_2\geq 0$.
Finally in Section \ref{ill-section} we  show that our results are sharp in the sense that the flow map of the
DMKP equation fail to be $C^2$ in $H^{s_1,0}(\rr^2)$) for $s_1 < -1/2$.

\section{Notations and Main Results}\label{2}
For the simplicity, throughout the paper we assume that $\beta=1$ (if $\beta\neq0$) and $\alpha=1$. Before stating our main result, we introduce our notations that are used in this paper.

We denote $\langle\cdot\rangle=1+|\cdot|$. The notation $A\lesssim B$ means that there exists the constant $C>0$ such that $A\leq CB$. Similarly, we will write $A\sim B$ to mean  $A\lesssim B$ and  $A\gtrsim B$.

For $n\in\N$, we denote by	$\what{f}$ the Fourier transform of $f$, defined as
\[
\what{f}(\omega)=\int_{\rr^n}f(x)\ee^{-\ii x\cdot \omega}\dd x.
\]
For $b,s_1,s_2\in\rr$, we denote $H^b=H^b(\rr)$, $\dot{H}^b=\dot{H}^b(\rr)$ and $H^{s_1,s_2}=H^{s_1,s_2}(\rr^2)$ as the nonhomogeneous Sobolev, the homogeneous Sobolev and the anisotropic Sobolev spaces, respectively, defined by
\begin{gather*}
H^b=\left\{f\in\mathscr{S}'(\rr);\;\|f\|_{H^b}=\|\langle\tau\rangle^b\what{f}(\tau)\|_{L^2_\tau}<\infty\right\},\\
\dot{H}^b=\left\{f\in\mathscr{S}'(\rr);\;\|f\|_{H^b}=\||\tau|^b\what{f}(\tau)\|_{L^2_\tau}<\infty\right\},\\
H^{s_1,s_2}=\left\{f\in\mathscr{S}'(\rr^2);\;\|f\|_{H^{s_1,s_2}}
=\|\langle\xi\rangle^{s_1}\langle\eta\rangle^{s_2}\what{f}(\xi,\eta)\|_{L^2_{\xi,\eta}}<\infty\right\}.
\end{gather*}
Let $U(\cdot)$ be the unitary group in $H^{s_1,s_2}$, $s_1, s_2\in\rr$, defining the free evolution of the KP equation \eqref{kp}, which is given by
\[
U(t)=\exp(\ii tP(D_x,D_y)),
\]
where $P(D_x,D_y )$ is the Fourier multiplier with symbol $P(\zeta)=P(\xi,\eta)=\xi^3-\vep\eta^2/\xi$, with $\vep=\pm1$.

We introduce a Bourgain space which is in relation with both the dissipative and dispersive parts of \eqref{dmkp} at the same time, we define this space by
\[
X^{b,s_1,s_2}=\left\{f\in\mathscr{S}'(\rr^3)\;:\;\|f\|_{X^{b,s_1,s_2}}<\infty\right\},
\]
equipped with the norm
\[
\|f\|_{X^{b,s_1,s_2}}=\left\|\langle\ii\sigma+\vr(\xi)\rangle^b\langle\xi\rangle^{s_1}\langle\eta\rangle^{s_2}
\what{f}(\tau,\zeta)\right\|_{L^2(\rr^3)},
\]
where $\sigma=\tau-P(\zeta)$ and $\vr(\xi)=\xi^4-\xi^2$.

We should note that $X^{b,s_1,s_2}$ is the intersection of the Bourgain space associated with the dispersive part of equation \eqref{dmkp} and Sobolev space. Indeed, one can easily see that
\[
\|f\|_{X^{b,s_1,s_2}}\approx
\|U(-t)f\|_{H^b_tH^{s_1,s_2}_{x,y}}+\|f\|_{L_t^2H^{s_1+2b,s_2}_{x,y}}.
\]
For $T >0$, we define the restricted spaces $X^{b,s_1,s_2}_T$ by the norm
\[
\|f\|_{X^{b,s_1,s_2}_T}=\inf_{f\in X^{b,s_1,s_2}}\left\{\|g\|_{X^{b,s_1,s_2}}\;:\:g(t)=f(t)\,\mbox{on}\,[0,T]\right\}.
\]
We denote by $W(\cdot)$ the semi-group associated with the free evolution of \eqref{dmkp},
\[
(W(t)f)^{\wedge_z}(\zeta)=\exp(\ii tP(\zeta)-t\vr(\xi)),\qquad f\in \mathscr{S}',\,z=(x,y),\,t\geq0.
\]
Also, we can extend $W$ to a linear operator defined on the whole real axis by setting
\[
(W(t)f)^{\wedge_z}(\zeta)=\exp(\ii tP(\zeta)-|t|\vr(\xi)),\qquad f\in \mathscr{S}',\,t\in\rr.
\]
By the Duhamel integral formulation, equation \eqref{dmkp} can be written
\begin{equation}\label{integral}
u(t)=W(t)\ff-\int_0^tW(t-t')\Lambda(u^2(t'))\dd t',\qquad t\geq0,
\end{equation}
where $\Lambda=\frac{1}{2}\partial_x+\partial_x^2$.

To prove the local existence result, we will apply a fixed point argument to a truncated version of \eqref{integral} which is defined on all the real axis by
\begin{equation}\label{trunc-int}
u(t)=\theta(t)W(t)\ff-\theta(t)\chi_{\rr^+}(t)\int_0^tW(t-t')\Lambda\left(\theta_T^2(t')u^2(t')\right)\dd t',\qquad t\geq0,
\end{equation}
where $t\in\rr$ and $\theta$ indicates a time cutoff function:
\[
\theta\in C_0^\infty(\rr),\quad {\rm supp}(\theta)\subset[-2,2],\quad\theta\equiv1\quad\mbox{on}\,\;[-1,1],
\]
and $\theta_T(\cdot)=\theta(\cdot/T)$.

We note that that if $u$ solves \eqref{trunc-int} then $u$ is a solution of \eqref{integral} on $[0,T ]$, $T\leq1$. Thus it is sufficient to solve \eqref{trunc-int} for a small time ($T \leq1$ is enough).

Let us now state our results.
\begin{theorem}\label{local}
Let $s_1>-1/2$, $s_2\geq0$, $s'_1\in(-1/2,\min\{0, s_1\}]$ and $\ff\in H^{s_1,s_2}$. Then there exist a time $T=T(\|\ff\|_{H^{s'_1,0}})>0$ and a unique solution $u$ of \eqref{dmkp} in
\begin{equation}
\mathscr{Y}_T=C([0,T],H^{s_1,s_2})\cap X_T^{1/2,s_1,s_2}.
\end{equation}
Moreover, $u$ belongs $C([0,T]; H^{s_1,s_2})\cap C((0,T]; H^{\infty,s_2})$ and the map $\ff\mapsto u$ is analytic from $H^{s_1,s_2}$ to $\mathscr{Y}_T$.
\end{theorem}
\begin{theorem}\label{ill-thm}
Let $s<-1/2$. Then it does not exist a time $T>0$ such that equation \eqref{dmkp} admits a unique solution in $C([0,T);H^{s,0})$ for any initial data in some ball of $H^{s,0}$ centered at the origin and such that the map $\ff\to u$ is $C^2$-differentiable at the origin from $H^{s,0}$ to $C([0,T ],H^{s,0})$.
\end{theorem}
\begin{theorem}\label{global}
Let $\beta=0$. Then Theorem \ref{local} holds for $s_1>-1$, $s_2\geq0$ and $\ff\in H^{s_1,s_2}$; and the corresponding local solution $u$ of the Cauchy problem \eqref{dmkp} extends globally in time.
\end{theorem}

\section{Linear Estimates}\label{3}
In this section we are going to obtain some appropriate linear estimates for \eqref{trunc-int}. The proofs of the linear estimates follow closely the proofs given in \cite{molinetribaud-1,molinetribaud,molinetribaud-2}.
In this section we study the linear operator $\theta V$.
\begin{lemma}\label{lem1}
Let $s_1,s_2\in\rr$, then for all $\ff\in H^{s_1,s_2}$, we have
\begin{equation}
\|\theta(t)W(t)\ff\|_{X^{1/2,s_1,s_2}}\lesssim\|\ff\|_{H^{s_1,s_2}}.
\end{equation}
\end{lemma}
\pf
By definition of $W(\cdot)$ and $X^{1/2,s_1,s_2}$, and by performing the change of variable $\tau\mapsto\sigma:=\tau-P(\zeta)$, we have
\begin{equation}\begin{split}
\|\theta(t)W(t)\ff\|_{X^{1/2,s_1,s_2}}&=
\left\|\left\langle\ii\sigma+\vr(\xi)\right\rangle^{1/2}
\langle|\xi|\rangle^{s_1}\langle|\eta|\rangle^{s_2}\left(\theta(t)\ee^{-|t|\vr(\xi)}
\what{\ff}(\zeta)\right)^{\wedge_t}(\tau)\right\|_{L^2(\rr^3)}\\
&=\left\|\langle|\xi|\rangle^{s_1}\langle|\eta|\rangle^{s_2}\what{\ff}(\zeta)\left\|\left\langle\ii \tau+\vr(\xi)\right\rangle^{1/2}
\left(\theta(t)\ee^{-|t|\vrx}\right)^{\wedge_t}(\tau)\right\|_{L^2_\tau}\right\|_{L^2_\zeta}\\
&\lesssim I+II,
\end{split}\end{equation}
where
\[
I=\left\|\langle|\xi|\rangle^{s_1}\langle|\eta|\rangle^{s_2}
\left\langle\vr(\xi)\right\rangle^{1/2}\what{\ff}(\zeta)
\|g_\xi(t)\|_{L^2_t}\right\|_{L^2_\zeta},
\]
\[
II=\left\|\langle|\xi|\rangle^{s_1}\langle|\eta|\rangle^{s_2}
\what{\ff}(\zeta)\|g_\xi(t)\|_{H_t^{1/2}}\right\|_{L^2_\zeta},
\]
and
\begin{equation}\label{gg}
g_\xi(t)=\theta(t)\ee^{-\vr(\xi)|t|}.
\end{equation}
\noindent\textbf{Contribution of $I$}.
When $|\xi|\geq\sqrt{2}$, we have $\vr(\xi)\geq2$, then we can obtain
\[
\|g_\xi\|_{L^2_t}\leq\|\ee^{-\vrx|t|}\|_{L^2_t}\sim|\vrx|^{-1/2}
\lesssim\frac{1}{\langle\vrx\rangle^{1/2}}.
\]
When $|\xi|\leq\sqrt{2}$, then $-1/4\leq\vrx\leq2$ implies that
\[
\|g_\xi\|_{L^2_t}\leq\left\|\ee^{|t|/2}\right\|_{L^2[-2,2]}\lesssim1\lesssim\langle\vrx\rangle^{-1/2}.
\]
Then we deduce that
\begin{equation}
I\lesssim\|\ff\|_{H^{s_1,s_2}}.
\end{equation}
\noindent\textbf{Contribution of $II$}.
When $|\xi|\geq\sqrt{2}$, we use the Young inequality to see that
\[\begin{split}
\|g_\xi\|_{H^{1/2}}&=\|\langle\tau\rangle^{1/2}\what{\theta}\ast(\ee^{-|t|\vrx})^{\wedge_t}(\tau)\|_{L^2_\tau}\\
&\lesssim\|\langle\tau\rangle^{1/2}\what{\theta}(\tau)\|_{L_\tau^1}\|\ee^{-|t|\vrx}\|_{L^2_t}
+\|\what{\theta}\|_{L^1_\tau}\|\ee^{-|t|\vrx)}\|_{\dot{H}_t^{1/2}}\\
&\lesssim\frac{1}{\langle\vrx\rangle^{1/2}}\lesssim1.
\end{split}\]
When $|\xi|\leq\sqrt{2}$, since  $|\vrx|\leq 2$, we have
\[
\|g_\xi\|_{H^{1/2}_t}\leq\sum_{j\geq0}\frac{2^j}{j!}\||t|^j\theta(t)\|_{H_t^{1/2}}\lesssim1.
\]
Since $
\||t|^j\theta(t)\|_{H_t^{1/2}}\leq\||t|^j\theta(t)\|_{H_t^1}\lesssim j,$ for $j\geq1$, therefore we deduce that
\begin{equation}
II\lesssim\|\ff\|_{H^{s_1,s_2}}.
\end{equation}
\fim
\begin{lemma}\label{lem2}
Let $0<\delta\leq1/2$ and $s_1,s_2\in\rr$, there exists $C=C_\delta>0$ such that for all $w\in  X^{-1/2+\delta,s_1-4\delta,s_2}$, we have
\begin{equation}\label{nonhomog}
\left\|\theta(t)\chi_{\rr^+}(t)\int_0^tW(t-t')w(t')\dd t'\right\|_{X^{1/2,s_1,s_2}}\leq C\|w\|_{X^{-1/2+\delta,s_1-4\delta,s_2}}.
\end{equation}
\end{lemma}
\pf
Let $b\in\rr$. For $\zeta\in\rr^2$ fixed, we define the following time-Sobolev space
\[
\mathscr{Y}_\zeta^b=\{w\in\mathscr{S}'(\rr^3);\;
\|w\|_{\mathscr{Y}_\zeta^b}=\|\langle\ii\tau+\vrx\rangle^b\what{w}(\tau,\zeta)\|_{L_\tau^2}<\infty\}.
\]
First we shall show that for $\zeta\in\rr^2$, $0<\delta\leq1/2$ and $w\in\mathscr{S}(\rr^3)$, the following estimate holds:
\begin{equation}\label{3.7-ami}
\|K_\zeta(t)\|_{\mathscr{Y}_\zeta^{1/2}}\lesssim\langle\xi\rangle^{-4\delta}\|w\|_{\mathscr{Y}_\zeta^{-1/2+\delta}},
\end{equation}
where
\[
K_\zeta(t)=\theta(t)\int_0^t\ee^{-|t-t'|\vrx}w(t',\zeta)\dd t'.
\]
By a simple calculation, similar to \cite{molinetribaud}, one can easily show that
\[
K_\zeta(t)=\theta(t)\int_\rr\frac{\ee^{\ii t\tau}-\ee^{-|t|\vrx}}{\ii\tau+\vrx}\what{w}(\tau,\zeta)\dd\tau.
\]
We split $K_\zeta$ into $K_\zeta=K_{1,0}+K_{1,\infty}+K_{2,0}+K_{2,\infty}$, where
\[
K_{1,0}=\theta(t)\int_{|\tau|\leq1}\frac{\ee^{\ii t\tau}-1}{\ii\tau+\vrx}\what{w}(\tau,\zeta)\dd\tau,\qquad
K_{1,\infty}=\theta(t)\int_{|\tau|\geq1}\frac{\ee^{\ii t\tau}}{\ii\tau+\vrx}\what{w}(\tau,\zeta)\dd\tau,
\]
\[
K_{2,0}=\theta(t)\int_{|\tau|\leq1}\frac{1-\ee^{-|t|\vrx}}{\ii\tau+\vrx}\what{w}(\tau,\zeta)\dd\tau,\qquad
K_{2,\infty}=\theta(t)\int_{|\tau|\geq1}\frac{\ee^{-|t|\vrx}}{\ii\tau+\vrx}\what{w}(\tau,\zeta)\dd\tau;
\]
and then we examine each $K_{\cdot,\cdot}$ in \eqref{3.7-ami}.

\noindent\textbf{Contribution of $K_{2,\infty}$}. In this case, since $|\tau|\geq1$, note that
\[
\|\langle\ii\tau+\vrx\rangle^{1/2}\what{K_{2,\infty}}\|_{L_\tau^2}
\leq
\left\|\langle\ii\tau+\vrx\rangle^{1/2}(g_\xi(t))^{\wedge_t}(\tau)\right\|_{L_\tau^2}
\left(\int_{|\tau|\geq1}\frac{\what{w}(\tau,\zeta)}{\langle\ii\tau+\vrx\rangle}\dd\tau\right),
\]
where $g_\xi$ is defined in \eqref{gg}. Exactly the same computations as in Lemma \ref{lem1} lead to
\[
\left\|\langle\ii\tau+\vrx\rangle^{1/2}(g_\xi(t))^{\wedge_t}(\tau)\right\|_{L_\tau^2}\lesssim1.
\]
Therefore, by the Cauchy-Schwarz inequality, we obtain
\[
\|K_{2,\infty}\|_{\mathscr{Y}_\zeta^{1/2}}
\lesssim
\left(\int_\rr\frac{|\what{w}(\tau,\zeta)|^2}{\langle\ii\tau+\vrx\rangle^{1-2\delta}}\dd\tau\right)^{1/2}
\left(\int_{|\tau|\geq1}\frac{\dd\tau}{\langle\ii\tau+\vrx\rangle^{1+2\delta}}\right)^{1/2}.
\]
When $|\xi|\geq\sqrt{2}$, a change of variable gives
\begin{equation}\label{k2inf}
\|K_{2,\infty}\|_{\mathscr{Y}_\zeta^{1/2}}
\lesssim
\langle\xi\rangle^{-4\delta}\|w\|_{\mathscr{Y}_\zeta^{-1/2+\delta}}.
\end{equation}
When $|\xi|\leq\sqrt{2}$, it follows $\langle\xi\rangle^{-4\delta}\sim1$; so that \eqref{k2inf} holds.

\noindent\textbf{Contribution of $K_{1,\infty}$}. In this case, by using the Young inequality, we see that
\[\begin{split}
\|K_{1,\infty}\|_{\mathscr{Y}_\zeta^{1/2}}&=\left\|\langle\ii\sigma+\vrx\rangle^{1/2}
\left|\what{\theta}(\tau')\ast\left(\frac{\what{w}(\tau',\zeta)}{|\ii\tau+\vrx|}\chi_{|\tau'|\geq1}\right)\right|(\tau)
\right\|_{L_\tau^2}\\
&\lesssim
\left\|\langle\tau'\rangle^{1/2}|\what{\theta}(\tau')|\ast
\left(\frac{\what{w}(\tau',\zeta)\chi_{|\tau'|\geq1}}{|\ii\tau'+\vrx|}\right)(\tau)\right\|_{L^2_\tau}
+\left\||\what{\theta}(\tau')|\ast
\left(\frac{\what{w}(\tau',\zeta)\chi_{|\tau'|\geq1}}{|\ii\tau'+\vrx|^{1/2}}\right)(\tau)\right\|_{L^2_\tau}\\
&\lesssim\|\langle\tau\rangle^{1/2}\what{\theta}(\tau)\|_{L^1_{\tau}}
\left\|\frac{\what{w}(\tau,\zeta)\chi_{|\tau|\geq1}}
{|\ii\tau'+\vrx|}\right\|_{L^2_{\tau}}
+\|\widehat{\theta}(\tau)\|_{L^1_{\tau}}\left\|\frac{\what{w}(\tau,\zeta)\chi_{|\tau|\geq1}}
{|\ii\tau+\vrx|^{1/2}}\right\|_{L^2_{\tau}}\\
&\lesssim
\left\|\frac{\what{w}(\tau,\zeta)}
{\langle\ii\tau+\vrx\rangle^{1/2}}\chi_{|\tau|\geq1}\right\|_{L^2_{\tau}}\\
&\lesssim\langle\xi\rangle^{-4\delta}\|w\|_{\mathscr{Y}_\zeta^{-1/2+\delta}}.
\end{split}\]
\noindent\textbf{Contribution of $K_{2,0}$}. First we notice that
\begin{equation}\label{k20-1}
\|K_{2,0}\|_{\mathscr{Y}_\zeta^{1/2}}\leq
\left(\int_{|\tau|\leq1}\frac{|\what{w}(\tau,\zeta)|}{|\ii\tau+\vrx|}\dd\tau\right)
\left\|\langle\ii\tau+\vrx\rangle^{1/2}\left(\theta(t)\left(1-\ee^{-|t|\vrx}\right)\right)^{\wedge_t}(\tau)\right\|_{L^2_\tau}.
\end{equation}
Now, as in the proof of Lemma \ref{lem1}, we consider two cases. When $|\xi|\geq\sqrt{2}$, we have $\vrx\geq2$, so that
\begin{equation}\label{k20-2}
\begin{split}
&\|\langle\ii\sigma+\vrx\rangle^{1/2}\left(\theta(t)\left(1-\ee^{-|t|\vrx}\right)\right)^{\wedge_t}(\tau)\|_{L_\tau^2}\\
&\qquad\lesssim \|\theta\|_{H^{1/2}_t}+\langle\vrx\rangle^{1/2}\|\theta\|_{L^2_t}+\|g_\zeta\|_{H_t^{1/2}}+
\langle\vrx\rangle^{1/2}\|g_\zeta\|_{L^2_t}\lesssim|\vrx|^{1/2}.
\end{split}\end{equation}
On the other hand, we have
\begin{equation}\label{k20-3}
\begin{split}
\int_{|\tau|\leq1}\frac{\langle\ii\tau+\vrx\rangle}{|\ii\tau+\vrx|^2}\dd\tau
&\approx
\int_{|\tau|\leq1}\frac{\dd\tau}{|\ii\tau+\vrx|^2}
+
\int_{|\tau|\leq1}\frac{\dd\tau}{|\ii\tau+\vrx|}\\
&\lesssim
\int_0^1\frac{\dd\tau}{\tau^2+\vr^2(\xi)}+\frac{1}{|\vrx|}\\
&\lesssim
\frac{1}{\vrx}\int_0^1\frac{1}{1+\left(\frac{\tau}{|\vrx|}\right)^2}\dd\left(\frac{\tau}{|\vrx|}\right)+\frac{1}{|\vrx|}
\lesssim\frac{1}{|\vrx|}.
\end{split}\end{equation}
From \eqref{k20-1}-\eqref{k20-3}, we deduce that
\[\begin{split}
\|K_{2,0}\|_{\mathscr{Y}_\zeta^{1/2}}&\lesssim|\vrx|^{1/2}
\left(\int_{|\tau|\leq1}\frac{|\what{w}(\tau,\zeta)|^2}{\langle\ii\tau+\vrx\rangle}\dd\tau\right)^{1/2}
\left(\int_{|\tau|\leq1}\frac{\langle\ii\tau+\vrx\rangle}{|\ii\tau+\vrx|^2}\dd\tau\right)^{1/2}\\
&\lesssim
\left(\int_{|\tau|\leq1}\frac{|\what{w}(\tau,\zeta)|^2}{\langle\ii\tau+\vrx\rangle}\dd\tau\right)^{1/2}
\lesssim
\langle\xi\rangle^{-4\delta}\|w\|_{\mathscr{Y}_\zeta^{-1/2+\delta}}.
\end{split}\]
When $|\xi|\leq\sqrt{2}$, then $|\vrx|\leq2$ and we have
\begin{equation}\label{k20-4}
\|\langle\ii\sigma+\vrx\rangle^{1/2}\left(\theta(t)\left(1-\ee^{-|t|\vrx}\right)\right)^{\wedge_t}(\tau)\|_{L_\tau^2}
\lesssim
\|\theta(t)\left(1-\ee^{-|t|\vrx}\right)\|_{H_t^{1/2}}.
\end{equation}
Then arguing again as in Lemma \ref{lem1}, we obtain that
\begin{equation}\label{k20-5}
\|\theta(t)\left(1-\ee^{-|t|\vrx}\right)\|_{H_t^{1/2}}
\leq\sum_{j\geq0}\frac{|\vrx|^j}{j!}\|t^j\theta(t)\|_{H_t^{1/2}}
\lesssim
|\vrx|\sum_{j\geq0}\frac{|\vrx|^j}{j!}\lesssim|\vrx|.
\end{equation}
From \eqref{k20-1} and \eqref{k20-3}-\eqref{k20-5}, we get
\[
\|K_{2,0}\|_{\mathscr{Y}_\zeta^{1/2}}\lesssim
\langle\xi\rangle^{-4\delta}\|w\|_{\mathscr{Y}_\zeta^{-1/2+\delta}}.
\]
\noindent\textbf{Contribution of $K_{1,0}$}. Since $K_{1,0}$ can be written as
\[
K_{1,0}=\theta(t)\sum_{j\geq1}\int_{|\tau|\leq1}\frac{(\ii t\tau)^j}{j!(\ii\tau+\vrx)}\what{w}(\tau,\zeta)\dd\tau,
\]
we deduce from the Cauchy-Schwarz inequality that
\[
\begin{split}
\|\langle\ii\tau+\vrx\rangle^{1/2}\what{K_{1,0}}(\tau)\|_{L_\tau^2}
&\lesssim
\sum_{j\geq1}\frac{1}{j!}\left(\|t^j\theta(t)\|_{H_t^{1/2}}+\langle\vrx\rangle^{1/2}\|t^j\theta(t)\|_{L_t^{1/2}}\right)
\int_{|\tau|\leq1}\frac{|\tau^j||\what{w}(\tau,\zeta)|}{|\ii\tau+\vrx|}\dd\tau\\
&\lesssim
\langle\vrx\rangle^{1/2}
\left(\int_\rr\frac{|\what{w}(\tau,\zeta)|^2}{\langle\ii\tau+\vrx\rangle}\dd\tau\right)^{1/2}
\left(\int_{|\tau|\leq1}\frac{|\tau^2|\langle\ii\tau+\vrx\rangle}{|\ii\tau+\vrx|^2}\dd\tau\right)^{1/2}\\
&\lesssim
\left(\int_\rr\frac{|\what{w}(\tau,\zeta)|^2}{\langle\ii\tau+\vrx\rangle}\dd\tau\right)^{1/2}.
\end{split}
\]
Finally, since for $\langle\ii\tau+\vrx\rangle^{1/2}\geq\langle\ii\tau+\vrx\rangle^{1-2\delta}\langle\vrx\rangle^{2\delta}$, we get
\[
\|K_{1,0}\|_{\mathscr{Y}_\zeta^{1/2}}\lesssim\langle\xi\rangle^{-4\delta}\|w\|_{\mathscr{Y}_\zeta^{-1/2+\delta}};
\]
which completes the proof of \eqref{3.7-ami}.

Now by definition of $X^{1/2,s_1,s_2}$, we see that
\[\begin{split}
&\left\|\theta(t)\chi_{\rr^+}(t)\int_0^tW(t-t')w(t')\dd t'\right\|_{X^{1/2,s_1,s_2}}\\
&=\left\|\langle\xi\rangle^{s_1}\langle\eta\rangle^{s_2}\langle\ii\sigma+\vrx\rangle^{1/2}
\left(\theta(t)\chi_{\rr^+}(t)\int_0^tW(t-t')w(t')\dd t'\right)^{\wedge_t}(\tau,\zeta)\right\|_{L^2(\rr^3)}.
\end{split}\]
We also note that
\[\begin{split}
&\left(\theta(t)\chi_{\rr^+}(t)\int_0^tW(t-t')w(t')\dd t'\right)^{\wedge_t}(\tau,\zeta)\\
&\qquad
=\left(\theta(t)\chi_{\rr^+}(t)\int_0^t\ee^{-|t-t'|\vrx}\ee^{\ii P(\zeta)(t-t')}\what{w}(t',\zeta)\dd t'\right)^{\wedge_t}(\tau,\zeta)\\
&\qquad
=\left(\theta(t)\chi_{\rr^+}(t)\int_0^t\ee^{-|t-t'|\vrx}\ee^{-\ii P(\zeta)t'}\left(U(t)w\right)^{\wedge_z}(t',\zeta)\dd t'\right)^{\wedge_t}(\tau,\zeta)\\
&\qquad
=\left(\theta(t)\chi_{\rr^+}(t)\int_0^t\ee^{-|t-t'|\vrx}\ee^{-\ii P(\zeta)t'}\left(w\right)^{\wedge_z}(t',\zeta)\dd t'\right)^{\wedge_t}(\tau-P(\zeta),\zeta);
\end{split}\]
and hence
\[\begin{split}
&\left\|\theta(t)\chi_{\rr^+}(t)\int_0^tW(t-t')w(t')\dd t'\right\|_{X^{1/2,s_1,s_2}}
\\&
\qquad
=\left\|\langle\xi\rangle^{s_1}\langle\eta\rangle^{s_2}
\left(\theta(t)\chi_{\rr^+}(t)\int_0^t\ee^{-|t-t'|\vrx}\left(U(-t)w\right)^{\wedge_z}(t',\zeta)\dd t'\right)^{\wedge_t}(\tau)\right\|_{L_\zeta^2(\mathscr{Y}^{1/2}_\zeta)}
\end{split}\]
Now define $v(t,\zeta)=(U(-t)w)^{\wedge_z}(t,\zeta)\in\mathscr{S}(\rr^3)$. Then by applying \eqref{3.7-ami}, we obtain
\[\begin{split}
&\left\|\theta(t)\chi_{\rr^+}(t)\int_0^tW(t-t')w(t')\dd t'\right\|_{X^{1/2,s_1,s_2}}
\\&
\qquad
\lesssim
\left\|\langle\xi\rangle^{s_1}\langle\eta\rangle^{s_2}
\|v\|_{\mathscr{Y}^{-1/2+\delta}_\zeta}\langle\xi\rangle^{-4\delta}\right\|_{L_\zeta^2}\\
&\qquad
\lesssim
\left\|\langle\xi\rangle^{s_1-4\delta}\langle\eta\rangle^{s_2}
\|\langle\ii\sigma+\vrx\rangle^{-1/2+\delta}v^{\wedge_t}(\tau)\|_{L_\tau^2}\right\|_{L_\zeta^2}\\
&\qquad
\lesssim
\left\|\langle\xi\rangle^{s_1-4\delta}\langle\eta\rangle^{s_2}\langle\ii\sigma+\vrx\rangle^{-1/2+\delta}
(U(-t)w)^{\wedge_{t,z}}(\tau,\zeta)\right\|_{L^2(\rr^3)}\\
&\qquad
\lesssim
\left\|\langle\xi\rangle^{s_1-4\delta}\langle\eta\rangle^{s_2}\langle\ii\sigma+\vrx\rangle^{-1/2+\delta}
\what{w}(\tau+P(\zeta),\zeta)\right\|_{L^2(\rr^3)}.
\end{split}\]
Finally, by performing a change of variable, we deduce \eqref{nonhomog}; and the proof of Lemma \ref{lem2} is complete.

\fim
\begin{lemma}\label{lem3}
Let $s_1,s_2\in\rr$ and $0<\delta\leq1/2$. Then for all $f\in X^{-1/2+\delta,s_1-4\delta,s_2}$, we have
\begin{equation}\label{lem3-in1}
N:t\longmapsto\int_0^tW(t-t')f(t')\dd t'\in C(\rr^+;H^{s_1,s_2}).
\end{equation}
Moreover, if $\{f_n\}$ is a sequence with $f_n\to0$ in $X^{-1/2+\delta,s_1-4\delta,s_2}$ as $n\to\infty$, then
\begin{equation}\label{lem3-in2}
\left\|\int_0^tW(t-t')f_n(t')\dd t'\right\|_{L^\infty(\rr^+;H^{s_1,s_2})}\longrightarrow 0.
\end{equation}
\end{lemma}
\pf
By Fubini theorem, and by the definition of $W(\cdot)$ we have
\begin{equation}\begin{split}
N(t)&=\int_0^t\left(\ee^{-|t-t'|\vrx}\ee^{\ii (t-t')P(\zeta)}(f(t'))^{\wedge_z}(\zeta)\right)^{\vee_\zeta}\dd t'\\
&=
U(-t)\left(\int_0^t\ee^{-|t-t'|\vrx}(g(t',\cdot))^{\wedge_z}(\zeta)\dd t'\right)^{\vee_\zeta},
\end{split}\end{equation}
where $g(t,z)=U(-t)f(t,\cdot)(z)$.
Since $U$ is  a strongly continuous unitary group in $L^2(\rr^2)$, it is enough to prove that
\[
F(\cdot,\zeta):t\in\rr^+\longmapsto
\langle\xi\rangle^{s_1}\langle\eta\rangle^{s_2}\int_0^t\ee^{-|t-t'|\vrx}(g(t',\cdot))^{\wedge_z}(\zeta)\dd t'
\]
is continuous from $\rr^+$ in $L^2_\zeta(\rr^2)$ for $f\in X^{-1/2+\delta,s_1-4\delta,s_2}$, $0<\delta\leq1/2$. We note that by Fubini theorem we have
\[
F(t,\zeta)=\langle\xi\rangle^{s_1}\langle\eta\rangle^{s_2}\int_\rr\what{g}(\tau,\zeta)\frac{\ee^{\ii t\tau}-\ee^{-t\vrx}}{\ii\tau+\vrx}\dd\tau.
\]
Fix $t_0\in\rr^+$ and define for all $t\in\rr$,
\[\begin{split}
H(t,\zeta):&=F(t,\zeta)-F(t_0,\zeta)\\
&=\langle\xi\rangle^{s_1}\langle\eta\rangle^{s_2}\int_\rr\frac{\what{g}(\tau,\zeta)}{\ii\tau+\vrx}
\left[\ee^{\ii t\tau}-\ee^{\ii t_0\tau}-\ee^{-t\vrx}+\ee^{-t_0\vrx}\right]\dd\tau.
\end{split}\]
We will use the Lebesgue dominated convergence theorem to show that
\begin{equation}\label{2.29}
\lim_{t\to t_0}\|H(t,\cdot)\|_{L^2(\rr^2)}=0.
\end{equation}
First we note that
\begin{equation}\label{2.32}
\lim_{t\to t_0}h(t,\tau,\zeta)=0,\qquad\mbox{a.e.}\,(\tau,\zeta)\in\rr^3,
\end{equation}
where
\begin{equation}\label{2.31}
h(t,\tau,\zeta)=\frac{\what{g}(\tau,\zeta)}{\ii\tau+\vrx}
\left[\ee^{\ii t\tau}-\ee^{\ii t_0\tau}-\ee^{-t\vrx}+\ee^{-t_0\vrx}\right].
\end{equation}
Moreover, since $t\to t_0$, we can suppose that $0\leq t\leq T$, and then,
\begin{equation}\label{2.33}
|h(t,\tau,\zeta)|\leq(2+\ee^{t/4}+\ee^{t_0/4})
\frac{|\what{g}(\tau,\zeta)|}{|\ii\tau+\vrx|}\lesssim
\frac{|\what{g}(\tau,\zeta)|}{|\ii\tau+\vrx|}.
\end{equation}
We deduce from the Cauchy-Schwarz inequality that
\[
\int_\rr\frac{|\what{g}(\tau,\zeta)|}{|\ii\tau+\vrx|}\dd\tau
\lesssim \left\|\frac{\langle\ii\tau+\vrx\rangle^{1/2-\delta}}{|\ii\tau+\vrx|}\right\|_{L^2_\tau}
\left\|\frac{\what{g}(\tau,\zeta)}{\langle\ii\tau+\vrx\rangle^{1/2-\delta}}\right\|_{L^2_\tau}.
\]
By the hypotheses on $g$, we deduce
\begin{equation}\label{2.34}
\int_\rr\frac{|\what{g}(\tau,\zeta)|}{|\ii\tau+\vrx|}\dd\tau\lesssim
\left\|\frac{\what{g}(\tau,\zeta)}{\langle\ii\tau+\vrx\rangle^{1/2-\delta}}\right\|_{L^2_\tau},
\end{equation}
for almost every $\zeta\in\rr^2$. We use \eqref{2.32}-\eqref{2.34} and the Lebesgue dominated
convergence theorem to conclude that
\begin{equation}\label{2.30}
\lim_{t\to t_0}H(\zeta,t)=0,\qquad \mathrm{a.e.}\quad\zeta\in\rr^2.
\end{equation}
Next we show that there exists $G\in L^2(\rr^2)$ such that
\begin{equation}\label{2.35}
|H(t,\zeta)|\leq|G(\zeta)|,
\end{equation}
for all $\zeta\in\rr^2$ and $t\in\rr^+$.

When $|\xi|\geq\sqrt{2}$, we get from the Cauchy-Schwarz inequality and  \eqref{2.33} that
\[
|H(t,\zeta)|\lesssim\langle\xi\rangle^{s_1}\langle\eta\rangle^{s_2}
\left\|\frac{\langle\ii\tau+\vrx\rangle^{1/2-\delta}}{|\ii\tau+\vrx|}\right\|_{L^2_\tau}
\left\|\frac{\what{g}(\tau,\zeta)}{\langle\ii\tau+\vrx\rangle^{1/2-\delta}}\right\|_{L^2_\tau}.
\]
Since $\vrx\geq2$, we have
\[
\left\|\frac{\langle\ii\tau+\vrx\rangle^{1/2-\delta}}{|\ii\tau+\vrx|}\right\|_{L^2_\tau}
\lesssim\left(\int_\rr\frac{1}{|\ii\tau+\vrx|^{1+2\delta}}\dd\tau\right)^{1/2}
\lesssim\langle\xi\rangle^{-4\delta},
\]
then using the hypotheses on $g$, we conclude that for all $t\in\rr^+$,
\[
|H(t,\zeta)|\lesssim\langle\xi\rangle^{s_1-4\delta}\langle\eta\rangle^{s_2}
\left\|\frac{\what{g}(\tau,\zeta)}{\langle\ii\tau+\vrx\rangle^{1/2-\delta}}\right\|_{L^2_\tau}
\in L^2(\rr^2),
\]
which proves \eqref{2.35} in this case. When $|\xi|\leq\sqrt{2}$, then we have $|\vrx|\leq2$, so that
\[\begin{split}
|H(t,\zeta)|&\lesssim\langle\xi\rangle^{s_1}\langle\eta\rangle^{s_2}
\int_\rr\frac{|\what{g}(\tau,\zeta)|}{|\ii\tau+\vrx|}
\left|\ee^{-t\vrx}-\ee^{-t_0\vrx}\right|\dd\tau+
\langle\xi\rangle^{s_1}\langle\eta\rangle^{s_2}
\int_\rr\frac{|\what{g}(\tau,\zeta)|}{|\ii\tau+\vrx|}
\left|\ee^{\ii t\tau}-\ee^{\ii t_0\tau}\right|\dd\tau\\
&=I+II.
\end{split}\]
We first evaluate $II$. Using the Cauchy-Schwarz inequality
\[\begin{split}
II&\leq|t-t_0|\langle\xi\rangle^{s_1}\langle\eta\rangle^{s_2}
\int_{|\tau|\leq1}\frac{|\tau||\what{g}(\tau,\zeta)|}{|\ii\tau+\vrx|}\dd\tau
+2\langle\xi\rangle^{s_1}\langle\eta\rangle^{s_2}
\int_{|\tau|\geq1}\frac{|\what{g}(\tau,\zeta)|}{|\ii\tau+\vrx|}\dd\tau\\
&\lesssim\langle\xi\rangle^{s_1-4\delta}\langle\eta\rangle^{s_2}
\left(\frac{|\what{g}(\tau,\zeta)|^2}{\langle\ii\tau+\vrx\rangle^{1-2\delta}}\dd\tau\right)^{1/2}
\left[\left(\int_{|\tau|\leq1}|\tau|^{1-2\delta}\dd\tau\right)^{1/2}+
\left(\int_{|\tau|\geq1}\langle\tau\rangle^{-1-2\delta}\dd\tau\right)^{1/2}\right]\\
&\lesssim\langle\xi\rangle^{s_1-4\delta}\langle\eta\rangle^{s_2}
\left(\frac{|\what{g}(\tau,\zeta)|^2}{\langle\ii\tau+\vrx\rangle^{1-2\delta}}\dd\tau\right)^{1/2}
\in L^2(\rr^2).
\end{split}\]
We next turn to $I$ and again use the Cauchy-Schwarz inequality to see that
\[
I\leq|t-t_0|\langle\xi\rangle^{s_1}\langle\eta\rangle^{s_2}
\left(\int_\rr\frac{|\what{g}(\tau,\zeta)|^2}{\langle\ii\tau+\vrx\rangle^{1-2\delta}}\dd\tau\right)^{1/2}
|\vrx|
\left(\int_\rr\frac{\langle\ii\tau+\vrx\rangle^{1-2\delta}}{|\ii\tau+\vrx|^2}\dd\tau\right)^{1/2},
\]
and we compute
\[\begin{split}
\left(\int_\rr\frac{\langle\ii\tau+\vrx\rangle^{1-2\delta}}
{|\ii\tau+\vrx|^2}\dd\tau\right)^{1/2}
&\lesssim
\left(\int_\rr\frac{1}{|\ii\tau+\vrx|^2}\dd\tau\right)^{1/2}
+
\left(\int_\rr\frac{1}{|\ii\tau+\vrx|^{1+2\delta}}\dd\tau\right)^{1/2}\\
&\lesssim\frac{1}{\sqrt{|\vrx|}}+\frac{1}{|\vrx|^\delta}.
\end{split}\]
Then, since $|\vrx|\leq2$, we conclude that
\[
I\lesssim \langle\xi\rangle^{s_1-4\delta}\langle\eta\rangle^{s_2} \left(\int_\rr\frac{|\what{g}(\tau,\zeta)|^2}{\langle\ii\tau+\vrx\rangle^{1-2\delta}}\dd\tau\right)^{1/2}
\in L^2(\rr^2).
\]
Thus \eqref{2.35} still remains true in this case. We use \eqref{2.30}, \eqref{2.35} and the dominated convergence theorem to prove \eqref{2.29}.

To show \eqref{lem3-in2} it suffices to  notice that one has
\[
\sup_{t\in\rr^+}\|F_n(t)\|_{L^2(\rr^2)}\lesssim\|f_n\|_{X^{-1/2+\delta,s_1-4\delta,s_2}},
\]
where $F_n$ is defined as $F$ with $g_n(t,z) =U(-t)f_n(t,\cdot)(z)$ instead of $g$. This completes the
proof.
\fim
\section{Bilinear Estimates}\label{4}
In this section, we are going to obtains suitable estimates for the nonlinear terms \eqref{dmkp}. Before stating this result,
we will give certain multilinear estimates which are necessary to treat the nonlinear term $\Lambda(u^2)$ in $X^{b,s_1,s_2}$.
\begin{lemma}[\cite{kojok,molinetribaud}]\label{lem4.3}
Let $u,v,w\in L^2(\rr^3)$ with compact support in $\{(x,y,t)\in\rr^3\;:\;|t|\leq T\}$. Then for $b>0$ and $c>0$ small enough there exists $\mu>0$ such that
\begin{equation}
\int_{\rr^6}\frac{|\what{u}(\tau,\zeta)||\what{v}(\tau_1,\zeta_1)||\what{w}(\tau_2,\zeta_2)|}
{\langle\sigma_1\rangle^{1/2}|\xi_1|^{3b+c}\langle\sigma_2\rangle^{1/2-b}}\dd\tau\dd\zeta\dd\tau_1\dd\zeta_1
\leq CT^\mu\|u\|_{L^2(\rr^3)}\|v\|_{L^2(\rr^3)}\|w\|_{L^2(\rr^3)},
\end{equation}
where
\begin{equation}\label{3zeta3}
\zeta=(\xi,\eta),\quad\zeta_1=(\xi_1,\eta_1),\quad\zeta_2=\zeta-\zeta_1
\end{equation}
and
\[
\sigma=\tau-P(\zeta),\quad \sigma_1=\tau_1-P(\zeta_1),\quad\sigma_2=\tau_2-P(\zeta_2).
\]
\end{lemma}
\begin{lemma}[\cite{kojok,molinetribaud}]\label{lem4.5}
Let $u,v,w\in L^2(\rr^3)$ with compact support in $\{(x,y,t)\in\rr^3\;:\;|t|\leq T\}$, $\epsilon>0$ and $a,b,c\in[0,1/2+\epsilon]$ such that $a+b+c\geq1+2\epsilon$. Then there exists $\mu>0$ such that
\begin{equation}
\int_{\rr^6}\frac{|\what{u}(\tau,\zeta)||\what{v}(\tau_1,\zeta_1)||\what{w}(\tau_2,\zeta_2)|}
{\langle\sigma\rangle^a\langle\sigma_1\rangle^b\langle\sigma_2\rangle^c}\dd\tau\dd\zeta\dd\tau_1\dd\zeta_1
\leq CT^\mu\|u\|_{L^2(\rr^3)}\|v\|_{L^2(\rr^3)}\|w\|_{L^2(\rr^3)}.
\end{equation}
\end{lemma}
\begin{theorem}\label{bil-theo}
Let $\delta>0$ small enough, $s_2\geq0$ and $s_1>-1/2$. For all $u,v\in X^{1/2,s_1,s_2}$ with compact support in time and included in the subset $\{(t,x,y);\;t\in[-T,T]\}$, there exists $\mu>0$ such that the following bilinear estimate holds
\begin{equation}
\|\Lambda(uv)\|_{X^{-1/2+\delta,s_1-4\delta,s_2}}\leq CT^\mu\|u\|_{X^{1/2,s_1,s_2}}\|v\|_{X^{1/2,s_1,s_2}}.
\end{equation}
\end{theorem}
\pf
We proceed by duality. It is equivalent to show that for $\delta >0$ small enough and for all $w\in X^{1/2-\delta,-s_1+4\delta,-s_2}$,
\begin{equation}\label{duality}
|\langle\Lambda(uv),w\rangle|\leq CT^\mu\|u\|_{X^{1/2,s_1,s_2}}\|v\|_{X^{1/2,s_1,s_2}}\|w\|_{X^{1/2-\delta,-s_1+4\delta,-s_2}}.
\end{equation}
Let $f$, $g$ and $h$ respectively defined by
\begin{gather}
\what{f}(\tau,\zeta)=\langle\ii(\tau-P(\zeta))+\vrx\rangle^{1/2}
\langle\xi\rangle^{s_1}\langle\eta\rangle^{s_2}\what{u}(\tau,\zeta),\\
\what{g}(\tau,\zeta)=\langle\ii(\tau-P(\zeta))+\vrx\rangle^{1/2}
\langle\xi\rangle^{s_1}\langle\eta\rangle^{s_2}\what{v}(\tau,\zeta),\\
\what{h}(\tau,\zeta)=\langle\ii(\tau-P(\zeta))+\vrx\rangle^{-1/2+\delta}
\langle\xi\rangle^{-s_1+4\delta}\langle\eta\rangle^{-s_2}\what{w}(\tau,\zeta).
\end{gather}
It is clear that
\[
\|u\|_{X^{1/2,s_1,s_2}}=\|f\|_{L^2(\rr^3)},\quad\|v\|_{X^{1/2,s_1,s_2}}=\|g\|_{L^2(\rr^3)}\quad\mbox{and}\quad
\|w\|_{X^{-1/2+\delta,-s_1+4\delta,-s_2}}=\|h\|_{L^2(\rr^3)}.
\]
Thus by Plancherel theorem, inequality \eqref{duality} is equivalent to
\begin{equation}\label{bilinear-main}
\begin{split}
\int_{\rr^6}&\frac{|q(\xi)||\what{f}(\tau_1,\zeta_1)||\what{g}(\tau_2,\zeta_2)||\what{h}(\tau,\zeta)|
\langle\xi\rangle^{s_1-4\delta}\langle\eta\rangle^{s_2}}
{\langle\ii\sigma+\vrx\rangle^{1/2-\delta}\langle\ii\sigma_1+\vr(\xi_1)\rangle^{1/2}
\langle\ii\sigma_2+\vr(\xi_2)\rangle^{1/2}\langle\xi_1\rangle^{s_1}\langle\xi_2\rangle^{s_1}
\langle\eta_1\rangle^{s_2}\langle\eta_2\rangle^{s_2}}\dd\tau\dd\tau_1\dd\zeta\dd\zeta_1\\
&\qquad\qquad\leq CT^\mu\|u\|_{L^2_{t,z}}\|v\|_{L^2_{t,z}}\|w\|_{L^2_{t,z}},
\end{split}\end{equation}
where $q(\xi)=|\xi|+\xi^2$. We can assume that $s_2=0$ and $s_1\leq0$, since in the case $s_1,s_2\geq0$, we have
\[
\frac{\langle\eta\rangle^{s_2}}{\langle\eta_1\rangle^{s_2}\langle\eta_2\rangle^{s_2}}\lesssim1\quad\mbox{and}\quad
\frac{\langle\xi\rangle^{s_1}}{\langle\xi_1\rangle^{s_1}\langle\xi_2\rangle^{s_1}}\lesssim1,
\]
for all $\xi_1,\xi,\eta_1,\eta\in\rr$. We note that it suffices to prove \eqref{bilinear-main} for $q(\xi)=\xi^2$.

Therefore setting $s=-s_1\geq0$, it is enough to estimate
\begin{equation}\label{bilinear-i}
I=\int_{\rr^6}
\frac{|\what{f}(\tau_1,\zeta_1)||\what{g}(\tau_2,\zeta_2)||\what{h}(\tau,\zeta)|\xi^2\langle\xi_1\rangle^s\langle\xi_2\rangle^s}
{\langle\ii\sigma+\vrx\rangle^{1/2-\delta}\langle\ii\sigma_1+\vr(\xi_1)\rangle^{1/2}
\langle\ii\sigma_2+\vr(\xi_2)\rangle^{1/2}
\langle\xi\rangle^{s+4\delta}}\dd\tau\dd\tau_1\dd\zeta\dd\zeta_1.
\end{equation}
By a symmetry argument we can restrict ourselves to the set
\[
A=\left\{(\tau_1,\zeta_1,\tau,\zeta)\in\rr^6;\;|\sigma_2|\leq|\sigma_1|\right\}.
\]
Let $\mathpzc{K}\gg4$. We divide $A$ into the following subregions:
\begin{gather*}
A_1=\{(\tau_1,\zeta_1,\tau,\zeta)\in A\;:\;|\xi|\leq\mathpzc{K},\;|\xi_1|\leq2\mathpzc{K}\},\\
A_2=\{(\tau_1,\zeta_1,\tau,\zeta)\in A\;:\;|\xi|\leq\mathpzc{K},\;|\xi_1|\geq2\mathpzc{K}\},\\
A_3=\{(\tau_1,\zeta_1,\tau,\zeta)\in A\;:\;|\xi|\geq\mathpzc{K},\;\min\{|\xi_1|,|\xi_2|\}\leq2\},\\
A_4=\{(\tau_1,\zeta_1,\tau,\zeta)\in A\;:\;|\xi|\geq\mathpzc{K},\;\min\{|\xi_1|,|\xi_2|\}\geq2\}.
\end{gather*}
\noindent\textbf{Case 1}.  Contribution of $A_1$ to $I$. In this case we have $|\xi_2|\lesssim1$ and we see that
\[
\frac{\xi^2\langle\xi_2\rangle^s\langle\xi_1\rangle^s}{\langle\xi\rangle^{s+4\delta}}\lesssim1;
\]
and hence,
\[\begin{split}
I&\lesssim\int_{\rr^6}
\frac{|\what{f}(\tau_1,\zeta_1)||\what{g}(\tau_2,\zeta_2)||\what{h}(\tau,\zeta)|}
{\langle\ii\sigma+\vrx\rangle^{1/2-\delta}\langle\ii\sigma_1+\vr(\xi_1)\rangle^{1/2}
\langle\ii\sigma_2+\vr(\xi_2)\rangle^{1/2}}\dd\tau\dd\tau_1\dd\zeta\dd\zeta_1\\
&\lesssim
\int_{\rr^6}
\frac{|\what{f}(\tau_1,\zeta_1)||\what{g}(\tau_2,\zeta_2)||\what{h}(\tau,\zeta)|}
{\langle\sigma\rangle^{1/2-\delta}\langle\sigma_1\rangle^{1/2}
\langle\sigma_2\rangle^{1/2}}\dd\tau\dd\tau_1\dd\zeta\dd\zeta_1.
\end{split}
\]
By applying Lemma \ref{lem4.5}, we deduce that
\[
I\lesssim T^\mu
\|u\|_{L^2_{t,z}}\|v\|_{L^2_{t,z}}\|w\|_{L^2_{t,z}}.
\]
\noindent\textbf{Case 2}.  Contribution of $A_2$ to $I$. Since we have, in this case, $|\xi|\leq|\xi_1|/2$, it follows that $|\xi_1|\sim|\xi-\xi_1|$. Therefore
\[\begin{split}
I&\lesssim\int_{\rr^6}
\frac{|\what{f}(\tau_1,\zeta_1)||\what{g}(\tau_2,\zeta_2)||\what{h}(\tau,\zeta)|\langle\xi_1\rangle^s\langle\xi_2\rangle^s}
{\langle\ii\sigma+\vrx\rangle^{1/2-\delta}\langle\ii\sigma_1+\vr(\xi_1)\rangle^{1/2}
\langle\ii\sigma_2+\vr(\xi_2)\rangle^{1/2}}\dd\tau\dd\tau_1\dd\zeta\dd\zeta_1\\
&\lesssim
\int_{\rr^6}
\frac{|\what{f}(\tau_1,\zeta_1)||\what{g}(\tau_2,\zeta_2)||\what{h}(\tau,\zeta)|}
{\langle\sigma\rangle^{1/2-\delta}\langle\sigma_1\rangle^{1/2-s/4}
\langle\sigma_2\rangle^{1/2-s/4}}\dd\tau\dd\tau_1\dd\zeta\dd\zeta_1.
\end{split}
\]
By applying again Lemma \ref{lem4.5}, for $s<1-2\delta$, we obtain that
\[
I\lesssim T^\mu
\|u\|_{L^2_{t,z}}\|v\|_{L^2_{t,z}}\|w\|_{L^2_{t,z}}.
\]
\noindent\textbf{Case 3}.  Contribution of $A_3$ to $I$. We first assume that $\min\{|\xi_1|,|\xi_2|\}=|\xi_1|$ and thus $2\leq|\xi_2|$ and $|\xi_2|\leq2+|\xi|\leq(2+C)|\xi|$, for $C>0$, and therefore $|\xi|\sim|\xi_2|$. It follows that
\[
\frac{\xi^2\langle\xi_1\rangle^s\langle\xi_2\rangle^s}{\langle\xi\rangle^{s+4\delta}}\lesssim|\xi|^{2-4\delta}.
\]
Since $\langle\ii\sigma+\vrx\rangle^{1/2-\delta}\gtrsim\langle\xi\rangle^{2-4\delta}$, it results that
\[\begin{split}
I&\lesssim\int_{\rr^6}
\frac{|\what{f}(\tau_1,\zeta_1)||\what{g}(\tau_2,\zeta_2)||\what{h}(\tau,\zeta)|}
{\langle\sigma_1\rangle^{1/2}
\langle\ii\sigma_2+\vr(\xi_2)\rangle^{1/2}}\dd\tau\dd\tau_1\dd\zeta\dd\zeta_1
\lesssim
\int_{\rr^6}
\frac{|\what{f}(\tau_1,\zeta_1)||\what{g}(\tau_2,\zeta_2)||\what{h}(\tau,\zeta)|}
{\langle\sigma_1\rangle^{1/2}
\langle\sigma_2\rangle^{1/2-\ell}|\xi_2|^{4\ell}}\dd\tau\dd\tau_1\dd\zeta\dd\zeta_1\\
&\lesssim
\int_{\rr^6}
\frac{|\what{f}(\tau_1,\zeta_1)||\what{g}(\tau_2,\zeta_2)||\what{h}(\tau,\zeta)|}
{\langle\sigma_1\rangle^{1/2}|\xi_1|^{4\ell}
\langle\sigma_2\rangle^{1/2-\ell}}\dd\tau\dd\tau_1\dd\zeta\dd\zeta_1,
\end{split}
\]
for any  $\ell\in(0,1/2)$. The estimate $I\lesssim T^\mu \|u\|_{L^2_{t,z}}\|v\|_{L^2_{t,z}}\|w\|_{L^2_{t,z}}$ is now derived from Lemma \ref{lem4.3}. The other case where $\min\{|\xi_1|, |\xi_2|\} =|\xi_2|$, follows exactly in the same manner.

\noindent\textbf{Case 4}.  Contribution of $A_4$ to $I$. In this case we need to divide $A_4$ in two regions defined by
\begin{gather*}
A_4^1=\{(\tau_1,\zeta_1,\tau,\zeta)\in A_4\;:\;|\xi|\geq\mathpzc{K}\min\{|\xi_1|,|\xi_2|\}\},\\
A_4^2=\{(\tau_1,\zeta_1,\tau,\zeta)\in A_4\;:\;|\xi|\leq\mathpzc{K}\min\{|\xi_1|,|\xi_2|\}\}.
\end{gather*}
\noindent\textbf{Case 4.1}.  Contribution of $A_4^1$ to $I$. By a symmetry argument we can assume that $\min\{|\xi_1|,|\xi_2|\}=|\xi_1|$. It follows $|\xi|\geq\mathpzc{K}|\xi_1|$. Thusly $|\xi_2|\leq|\xi_1|+|\xi|\lesssim|\xi|$ and $|\xi|\leq|\xi_1|+|\xi_2|\leq|\xi|/\mathpzc{K}+|\xi_2|$ and consequently, $|\xi|\sim|\xi_2|$. It results that
\[
\frac{\xi^2\langle\xi_1\rangle^s\langle\xi_2\rangle^s}{\langle\xi\rangle^{s+4\delta}}
\lesssim|\xi|^{2-4\delta}|\xi_1|^s.
\]
Hence $\langle\ii\sigma_2+\vr(\xi_2)\rangle^{1/2}\gtrsim\langle\sigma_2\rangle^{1/2-s/4-\ell}|\xi_2|^{4\ell}
\gtrsim\langle\sigma_2\rangle^{1/2-s/4-\ell}|\xi_1|^{4\ell}$ gives
\[\begin{split}
I&\lesssim\int_{\rr^6}
\frac{|\what{f}(\tau_1,\zeta_1)||\what{g}(\tau_2,\zeta_2)||\what{h}(\tau,\zeta)|\;|\xi_1|^s}
{\langle\sigma_1\rangle^{1/2}
\langle\ii\sigma_2+\vr(\xi_2)\rangle^{1/2}}\dd\tau\dd\tau_1\dd\zeta\dd\zeta_1
\lesssim
\int_{\rr^6}
\frac{|\what{f}(\tau_1,\zeta_1)||\what{g}(\tau_2,\zeta_2)||\what{h}(\tau,\zeta)|}
{\langle\sigma_1\rangle^{1/2}
\langle\sigma_2\rangle^{1/2-s/4-\ell}|\xi_2|^{4\ell}}\dd\tau\dd\tau_1\dd\zeta\dd\zeta_1\\
&\lesssim
\int_{\rr^6}
\frac{|\what{f}(\tau_1,\zeta_1)||\what{g}(\tau_2,\zeta_2)||\what{h}(\tau,\zeta)|}
{\langle\sigma_1\rangle^{1/2}|\xi_1|^{4\ell}
\langle\sigma_2\rangle^{1/2-s/4-\ell}}\dd\tau\dd\tau_1\dd\zeta\dd\zeta_1,
\end{split}
\]
for any  $3s/4<\ell<1/2-s/4$. Therefore a use of Lemma \ref{lem4.3} provides a good bound for $I$ in this case.
\noindent\textbf{Case 4.2}.  Contribution of $A_4^2$ to $I$. To estimate $I$ in this case we need to split $A_4^1$ into the following two subregions:
\begin{gather*}
A_4^{21}=\{(\tau_1,\zeta_1,\tau,\zeta)\in A_4^1\;:\;|\sigma_1|\geq|\sigma|\},\\
A_4^{22}=\{(\tau_1,\zeta_1,\tau,\zeta)\in A_4^1\;:\;|\sigma|\geq|\sigma_1|\}.
\end{gather*}
\noindent\textbf{Case 4.21}.  Contribution of $A_4^{21}$ to $I$. In this case, by a symmetry argument we assume that $\min\{|\xi_1|,|\xi_2|\}=|\xi_1|$. We have $|\xi|\lesssim|\xi_1|$, $|\xi|\lesssim|\xi_2|$ and $\langle\ii\sigma_1+\vr(\xi_1)\rangle\gtrsim\langle\ii\sigma+\vrx\rangle$, and therefore we obtain
\[\begin{split}
I&\lesssim\int_{\rr^6}
\frac{|\what{f}(\tau_1,\zeta_1)||\what{g}(\tau_2,\zeta_2)||\what{h}(\tau,\zeta)|
\langle\xi\rangle^{2-s-4\delta}\langle\xi_1\rangle^s\langle\xi_2\rangle^s}
{\langle\sigma\rangle^{1/2}\langle\ii\sigma_1+\vr(\xi_1)\rangle^{1/2-\delta}
\langle\ii\sigma_2+\vr(\xi_2)\rangle^{1/2}}\dd\tau\dd\tau_1\dd\zeta\dd\zeta_1\\
&\lesssim
\int_{\rr^6}
\frac{|\what{f}(\tau_1,\zeta_1)||\what{g}(\tau_2,\zeta_2)||\what{h}(\tau,\zeta)|
\langle\xi_1\rangle^{2-4\delta}\langle\xi_2\rangle^s}
{\langle\sigma\rangle^{1/2}\langle\xi_1^4\rangle^{1/2-\delta}
\langle\sigma_2\rangle^{1/2-s/4-\ell}\langle\xi_2^4\rangle^{s/4+\ell}}\dd\tau\dd\tau_1\dd\zeta\dd\zeta_1\\
&\lesssim
\int_{\rr^6}
\frac{|\what{f}(\tau_1,\zeta_1)||\what{g}(\tau_2,\zeta_2)||\what{h}(\tau,\zeta)|}
{\langle\sigma\rangle^{1/2}
\langle\sigma_2\rangle^{1/2-s/4-\ell}|\xi_2|^{4\ell}}\dd\tau\dd\tau_1\dd\zeta\dd\zeta_1\\
&\lesssim
\int_{\rr^6}
\frac{|\what{f}(\tau_1,\zeta_1)||\what{g}(\tau_2,\zeta_2)||\what{h}(\tau,\zeta)|}
{\langle\sigma\rangle^{1/2}
\langle\sigma_2\rangle^{1/2-s/4-\ell}|\xi|^{4\ell}}\dd\tau\dd\tau_1\dd\zeta\dd\zeta_1,
\end{split}
\]
for any $3s/4<\ell<1/2-s/4$, and thus we can apply Lemma \ref{lem4.3} to estimate $I$ in this case.

\noindent\textbf{Case 4.22}.  Contribution of $A_4^{22}$ to $I$.
In this case, by a symmetry argument we again assume that $\min\{|\xi_1|,|\xi_2|\}=|\xi_1|$. We have $|\xi|\lesssim|\xi_1|$ and  $|\xi|\lesssim|\xi_2|$, and therefore it follows that
\[\begin{split}
I&\lesssim\int_{\rr^6}
\frac{|\what{f}(\tau_1,\zeta_1)||\what{g}(\tau_2,\zeta_2)||\what{h}(\tau,\zeta)|
\langle\xi\rangle^{2-s-4\delta}\langle\xi_1\rangle^s\langle\xi_2\rangle^s}
{\langle\ii\sigma+\vrx\rangle^{1/2-\delta-s/4}\langle\sigma\rangle^{s/4}
\langle\ii\sigma_1+\vr(\xi_1)\rangle^{1/2}
\langle\ii\sigma_2+\vr(\xi_2)\rangle^{1/2}}\dd\tau\dd\tau_1\dd\zeta\dd\zeta_1\\
&\lesssim
\int_{\rr^6}
\frac{|\what{f}(\tau_1,\zeta_1)||\what{g}(\tau_2,\zeta_2)||\what{h}(\tau,\zeta)|
\langle\xi_1\rangle^{2-s-4\delta}\langle\xi_1\rangle^s\langle\xi_2\rangle^s}
{\langle\xi^4\rangle^{1/2-\delta-s/4}\langle\sigma\rangle^{s/4}
\langle\sigma_1\rangle^{1/2-s/4}\langle\xi_1^4\rangle^{s/4}
\langle\sigma_2\rangle^{1/2-s/4}\langle\xi_2^4\rangle^{s/4}}\dd\tau\dd\tau_1\dd\zeta\dd\zeta_1\\
&\lesssim
\int_{\rr^6}
\frac{|\what{f}(\tau_1,\zeta_1)||\what{g}(\tau_2,\zeta_2)||\what{h}(\tau,\zeta)|}
{\langle\sigma\rangle^{s/4}\langle\sigma_1\rangle^{1/2-s/4}\langle\sigma_2\rangle^{1/2-s/4-\ell}|\xi_2|^{4\ell}}
\dd\tau\dd\tau_1\dd\zeta\dd\zeta_1\\
&\lesssim
\int_{\rr^6}
\frac{|\what{f}(\tau_1,\zeta_1)||\what{g}(\tau_2,\zeta_2)||\what{h}(\tau,\zeta)|}
{\langle\sigma_1\rangle^{1/2}
\langle\sigma_2\rangle^{1/2-s/4-\ell}|\xi_1|^{4\ell}}\dd\tau\dd\tau_1\dd\zeta\dd\zeta_1,
\end{split}
\]
for any $\ell\in(3s/4,1/2-s/4)$. Finally by applying Lemma \ref{lem4.3}, we estimate $I$ in this case.

This completes the proof of Theorem \ref{bil-theo}.
\fim
\begin{corollary}\label{bil-1-cor}
Let $\delta>$ sufficiently small, $s_1\geq s_1'>-1/2$ and $s_2\geq0$. Then for $u,v\in X^{1/2,s_1,s_2}$, compact supported in time in $\{(t,x,y)\in\rr^3\;:\;t\in[-T,T]\}$, there exists $\mu>0$ such that
\[\begin{split}
\|\Lambda(uv)\|_{X^{-1/2+4\delta,s_1-4\delta,s_2}}
&\lesssim
T^\mu
\left(
\|u\|_{X^{1/2,s_1',0}}\|v\|_{X^{1/2,s_1,s_2}}+
\|u\|_{X^{1/2,s_1,s_2}}\|v\|_{X^{1/2,s_1',0}}\right.\\&\left.\quad+
\|u\|_{X^{1/2,s_1,0}}\|v\|_{X^{1/2,s_1',s_2}}+
\|u\|_{X^{1/2,s_1',s_2}}\|v\|_{X^{1/2,s_1,0}}
\right).
\end{split}\]
\end{corollary}
\pf
The proof is a direct consequence of Theorem \ref{bil-theo} together with the following inequalities:
\begin{gather*}
\langle\xi\rangle^{s_1}\leq\langle\xi\rangle^{s_1'}\langle\xi_1\rangle^{s_1-s_1'}
+\langle\xi\rangle^{s_1'}\langle\xi-\xi_1\rangle^{s_1-s_1'},\quad s_1\geq s_1',\\
\langle\eta\rangle^{s_2}\leq\langle\eta_1\rangle^{s_2}+\langle\eta-\eta_1\rangle^{s_2}.
\end{gather*}
\fim
\begin{remark}\label{remark-g}
 It is noteworthy that by an argument similar to Theorem \ref{bil-theo}, one can show that  the bilinear estimate of Theorem  \ref{bil-theo} holds for $s_1>-1$ and $s_2\geq0$, if $\beta=0$ in \eqref{dmkp}.
\end{remark}
\section{Existence}\label{5}
Now we are ready to prove Theorem \ref{local} and Theorem \ref{global}.

\noindent\textbf{Proof of Theorem \ref{local}.}\quad
Let $\ff\in H^{s_1,s_2}$ with s$_1 >-1/2$, $s_2\geq0$ and $s_1'\in(-1/2, \min\{0, s_1\}]$. We suppose that $T\leq1$,
if $u$ is a solution of the integral equation $u = \Phi(u)$ with
\begin{equation}\label{operator}
\Phi(u)=\theta(t)\left(W(t)\ff-\chi_{\rr^+}(t)\int_0^tW(t-t')\Lambda\left(\theta_T^2(t')u^2(t')\right)\dd t'\right).
\end{equation}
then $u$ solve the DMKP equation on $[0,T/2]$. We introduce the Bourgain spaces defined by
\begin{gather}
\mathscr{Z}_1=\left\{u\in X^{1/2,s_1,s_2}\;:\;\|u\|_{\mathscr{Z}_1}=\|u\|_{ X^{1/2,s_1,0}}+\kappa_1\|u\|_{ X^{1/2,s_1,s_2}}<\infty\right\},\\
\mathscr{Z}_2=\left\{u\in X^{1/2,s_1,0}\;:\;\|u\|_{\mathscr{Z}_2}=\|u\|_{ X^{1/2,s_1',0}}+\kappa_2\|u\|_{ X^{1/2,s_1,0}}<\infty\right\},
\end{gather}
where
\[
\kappa_1=\frac{\|\ff\|_{H^{s_1,0}}}{\|\ff\|_{H^{s_1,s_2}}},\qquad
\kappa_2=\frac{\|\ff\|_{H^{s_1',0}}}{\|\ff\|_{H^{s_1,0}}}.
\]
We show that there exist $T_1 =T_1(H^{s_1,0})$ and a solution $u$ of \eqref{operator} in a ball of $\mathscr{Z}_1$, and then we solve \eqref{operator} in $\mathscr{Z}_2$ in order to check that the time of existence $T =T (H^{s_1',0})$.

First, by Lemmas \ref{lem1} and \ref{lem2}, we have
\begin{gather*}
\|\Phi (u)\|_{X^{1/2,s_1,0}}\lesssim
\|\ff\|_{H^{s_1,0}}+
T^{2\mu}\|\Lambda(\theta_T^2u^2)\|_{X^{-1/2+\delta,s_1-4\delta,0}},\\
\|\Phi (u)\|_{X^{1/2,s_1,s_2}}\lesssim
\|\ff\|_{H^{s_1,s_2}}+
T^{2\mu}\|\Lambda(\theta_T^2u^2)\|_{X^{-1/2+\delta,s_1-4\delta,s_2}}.
\end{gather*}
By Theorem  \ref{bil-theo}, Corollary \ref{bil-1-cor}, Leibniz rule for fractional derivative and Sobolev inequalities in time, we can deduce
\begin{gather*}
\|\Phi (u)\|_{X^{1/2,s_1,0}}\lesssim
\|\ff\|_{H^{s_1,0}}+
T^\mu\|u\|_{X^{1/2,s_1,0}}^2,\\
\|\Phi (u)\|_{X^{1/2,s_1,s_2}}\lesssim
\|\ff\|_{H^{s_1,s_2}}+
T^\mu\|u\|_{X^{1/2,s_1,s_2}}\|u\|_{X^{1/2,s_1,0}};
\end{gather*}
and consequently we obtain
\begin{equation}\label{exist-1}
\|\Phi(u)\|_{\mathscr{Z}_1}\leq
C\left(\|\ff\|_{H^{s_1,0}}+\kappa_1\|\ff\|_{H^{s_1,s_2}}\right)+CT^\mu\|u\|_{\mathscr{Z}_1}^2.
\end{equation}
Analogously, we can get
\begin{equation}\label{exist-2}
\|\Phi(u)-\Phi(v)\|_{\mathscr{Z}_1}\leq CT^\mu
\|u-v\|_{\mathscr{Z}_1}\|u+v\|_{\mathscr{Z}_1}.
\end{equation}
Hence by setting
\begin{equation}
T_1=\left[4C^2(\|\ff\|_{H^{s_1,0}}+\kappa_1\|\ff\|_{H^{s_1,s_2}})\right]^{-2/\mu}
=
\left[8C^2\|\ff\|_{H^{s_1,0}}\right]^{-2\mu},
\end{equation}
we can deduce from \eqref{exist-1} and \eqref{exist-2} that $\Phi$ is strictly contractive on the ball of radius
\[
2C(\|\ff\|_{H^{s_1,0}}+\kappa_1\|\ff\|_{H^{s_1,s_2}})
\]
in $\mathscr{Z}_1$ This proves the existence of a unique solution $u_1$ to \eqref{operator} in $X^{1/2,s_1,s_2}$ with $T_1$ defined above. On the other hand, Since $\ff\in H^{s_1,s_2}$ , it follows that $\theta(\cdot)W(\cdot)\ff\in C([0,T_1],H^{s_1,s_2})$, moreover since $u_1\in X^{1/2,s_1,s_2}$ , we can deduce from Theorem \ref{bil-theo} that $\Lambda(u_1^2)\in X^{-1/2+\delta,s_1-4\delta,s_2}$ and from Lemma \ref{lem3},  we obtain that
\begin{equation}
t\longmapsto\int_0^tW(t-t')\Lambda(u_1^2)\dd t'\in C([0,T_1];H^{s_1,s_2}).
\end{equation}
Thus $u_1$ belongs $C([0,T_1],H^{s_1,s_2})$.

An argument as above in $\mathscr{Z}_2$ shows that $\Phi$ is also strictly contractive on the ball of radius
\[
2C(\|\ff\|_{H^{s_1',0}}+\kappa_2\|\ff\|_{H^{s_1,0}})
\]
in $\mathscr{Z}_2$ with
\[
T_2=\left[4C^2(\|\ff\|_{H^{s_1',0}}+\kappa_2\|\ff\|_{H^{s_1,0}})\right]^{-1/\mu}.
\]
Therefore by definition of $\kappa_2$, it follows that $T_2=T_2(\|\ff\|_{H^{s_1',0}})$; which it follows that there exists a unique solution $u_1$ of \eqref{operator} in $C([0,T_2];H^{s_1,0})\cap X^{1/2,s_1,0}$. If we indicate
by $T^\ast= T_{\max}$ the maximum time of the existence in $\mathscr{Z}_1$ then by uniqueness, we have $u_1 = u_2$ on
$[0, \min\{T_2,T^\ast\})$ and this gives that $T^\ast\geq T_2(\|\ff\|_{H^{s_1',0}})$.

The continuity of map $\ff\mapsto u$ from $H^{s_1,s_2}$ to $X^{1/2,s_1,s_2}$ follows from classical argument, while
the continuity from $H^{s_1,s_2}$ to $C([0,T_1],H^{s_1,s_2} )$ follows again from Lemma \ref{lem3}. The analyticity
of the flow-map is a direct consequence of the implicit function theorem.

The uniqueness of the solution to the truncated integral equation \eqref{operator} is consequence of the contraction argument. We deduce  the uniqueness of the solution to the integral equation \eqref{integral} by using the ideas of \cite{molinetribaud-2}.

Let $u,v\in X^{1/2,s_1,s_2}$ be two solutions of the integral equation \eqref{integral} on the time interval $[0,T]$ and let $\ti{u}-\ti{v}$ be an extension of $u-v$ in $X^{1/2,s_1,s_2}$ such that
\[
\|\ti{u}-\ti{v}\|_{X^{1/2,s_1,s_2}}\lesssim\|u-v\|_{X^{1/2,s_1,s_2}_\kappa}
\]
with $0\leq\kappa\leq T/2$. It results by Lemmas \ref{lem1} and \ref{lem2} that
\[\begin{split}
\|u-v\|_{X^{1/2,s_1,s_2}_\kappa}
&\leq
\left\|\theta(t)\chi_{\rr^+}(t)\int_0^tW(t-t')\Lambda(\theta_\kappa^2(t')(\ti{u}^2-\ti{v}^2)(t'))\dd t'\right\|_{X^{1/2,s_1,s_2}}\\
&\leq\left\|\Lambda(\theta_\kappa^2(t')(\ti{u}^2-\ti{v}^2)(t'))\right\|_{X^{-1/2+\delta,s_1-4\delta,s_2}}\\
&\leq C\kappa^{\mu/2}
\|u+v\|_{X^{1/2,s_1,s_2}}\|\ti{u}-\ti{v}\|_{X^{1/2,s_1,s_2}}\\
&\leq2C\kappa^{\mu/2}\left(\|u\|_{X^{1/2,s_1,s_2}}+\|v\|_{X^{1/2,s_1,s_2}}\right)\|u-v\|_{X_\kappa^{1/2,s_1,s_2}}.
\end{split}\]
for some $\mu>0$. By considering $\kappa\leq\left[4C(\|u\|_{X_T^{1/2,s_1,s_2}}+\|v\|_{X_T^{1/2,s_1,s_2}})\right]^{-\mu/2}$, it follows that  $u\equiv v$ on $[0,\kappa]$. Iterating this argument, one extends the uniqueness result on the whole time interval $[0,T]$.
\fim

A proof of Theorem \ref{global} is now in sight.

\noindent\textbf{Proof of Theorem \ref{global}.}\quad
The local existence is obtained by an argument similar to Theorem \ref{local} and Remark \ref{remark-g}. To show the global existence when $\beta=0$, we note that $\partial_x(u^2)\in X^{-1/2+4\delta,s_1-\delta,s_2}$. Therefore by Lemma \ref{lem3}, we obtain that
\begin{equation}
t\longmapsto\int_0^t W(t-t')\partial_x(u^2(t'))\dd t'\in C([0,T]; H^{s_1+\epsilon,s_2}).
\end{equation}
Note that $$W(\cdot)\ff\in C([0,+\infty;H^{s_1,s_2})\cap C((0,+\infty);H^{\infty,s_2});$$ and consequently
$$u\in C([0,T];H^{s_1,s_2})\cap C((0,T);H^{s_1+\epsilon,s_2}).$$
Noting that $T = T (\|\ff\|_{H^{s_1',0}})$ with $s_1'>-1$ and using the uniqueness result, we deduce by induction that $u \in C((0,T];H^{\infty,s_2})$. This allows us to take the $L^2$-scalar product of the DMKP equation with $u$, which shows that $t\mapsto\|u(t)\|_{L^2}$ is nonincreasing on $(0,T]$. Since the time of
local existence $T$ only depends on  $\|\ff\|_{H^{s_1',0}}$, this clearly gives that the solution is global in time.
\fim
\section{Ill-poseness}\label{ill-section}
In this section, we prove the ill-posedness result for the DMKP equation stated in Theorem \ref{ill-thm}. We start by constructing a sequence of initial data $\{\ff_n\}_n$ which will ensure the nonregularity of the map $\ff\to u$ from $H^{s,0}$ to $C([0,T ],H^{s,0})$ for$ s <-1/2$.

\noindent\textbf{Proof of Theorem \ref{ill-thm}.}\quad
Let $u$ be a solution of \eqref{dmkp}. Then we have
\[
u(x,y,t,\ff)=W(t)\ff(x,y)-\int_0^tW(t-t')\Lambda(u^2(x,y,t',\ff))\dd t'.
\]
We will argue by contradiction and suppose that the map $\ff\to u$ is $C^2$. Since $u(x,y,0,\ff)=0$, it is straightforward to verify  that
\[
u_1(x,y,t)=\frac{\partial u}{\partial\ff}(x,y,t,0)[h]=W(t)h,
\]
\[
u_2(x,y,t)=\frac{\partial^2 u}{\partial\ff^2}(x,y,t,0)[h,h]=-\int_0^tW(t-t')\Lambda\left(\left(W(t')h\right)^2\right)\dd t'.
\]
The assumption of $C^2$-regularity of the map solution implies that
\begin{equation}
\|u_j(\cdot,\cdot,t)\|_{H^{s,0}}\lesssim
\|h\|_{H^{s,0}}^j,\qquad j=1,2,\quad\forall h\in H^{s,0}.
\end{equation}
First recall the definitions of $\zeta_1$, $\zeta$ and $\zeta_2$ in \eqref{3zeta3}. A straightforward calculation reveals that
\[
(u_2(\cdot,\cdot,t))^{\wedge_z}(\zeta)=
\ii(\xi+\xi^2)\ee^{\ii tP(\zeta)}\int_{\rr^2}\what{\ff}(\zeta_1)\what{\ff}(\zeta_2)\frac{\ee^{-t(\vr(\xi_1)+\vr(\xi_2))}\ee^{\ii t\cchi}-\ee^{-t\vrx}}{\M+\ii\cchi}\dd\zeta_1,
\]
where $\cchi=P(\zeta_1)+P(\zeta_2)-P(\zeta)$ and $\M=\vr(\xi_1)+\vr(\xi_2)-\vr(\xi)$. Note that from definitions of $P(\zeta)$ and $\vrx$, it is readily seen that
\[
\cchi=\mathscr{R}(\xi,\eta,\xi_1,\eta_1)=3\xi\xi_1\xi_2+\vep\frac{(\eta\xi_1-\eta_1\xi)^2}{\xi\xi_1\xi_2}
\]
and
\[
\M=-2\xi_1\xi_2(\xi_1^2-\xi\xi_1+2\xi^2-1).
\]
We choose now a sequence of initial data $\{\ff_N\}_N$, $N>0$, defined through its Fourier transform
by
\[
\what{\ff_N}(\xi,\eta)=N^{-3/2-s}(\chi_{A_N}(\xi,\eta)+\chi_{B_N}(\xi,\eta))
\]
where $N$ is a positive parameter such that $N\gg 1$, and $A_N$ , $B_N$ are the rectangles in $\rr^2$ defined by
\[
A_N=[N/2,N]\times\left[-6N^2,6N^2\right],\qquad B_N=[N,2N]\times\left[2N^2,3N^2\right].
\]
Note first that $\|\ff_N\|_{H^{s,0}}\sim1$. Let us denote by $u_{2,N}$ the sequence of the second iteration
$u_2$ associated with $\ff_N$. Hence it is readily seen that
\begin{equation}\label{norm-hs}
\|u_{2,N}\|_{H^{s,0}}^2\gtrsim
N^{-4s-6}\int_{\rr^2}(|\xi|+\xi^2)^2(1+|\xi|^2)^s\left|\int_{k_\zeta}\K\dd\zeta_1\right|^2\dd\zeta,
\end{equation}
where
\[
\K=\frac{\ee^{\vr(\xi_1)+\vr(\xi_2)}\ee^{\ii t\cchi}-\ee^{-t\vrx}}{\M+\ii\cchi}
\]
and
\[
k_\zeta=\{\zeta_1\;:\;\zeta_1\in B_N,\;\zeta_2\in A_N\}\cup\{\zeta_1\;:\;\zeta_1\in A_N,\;\zeta_2\in B_N\}.
\]
Now the definition of $\mathscr{M}$ shows that $|\M|\lesssim N^4$. On the other hand,  by Lemma 7.1 in \cite{kojok},  we deduce from the inequality
\[
|\cchi|\leq3(1+\varepsilon)|\xi\xi_1\xi_2|+\left|3\xi\xi_1\xi_2-\frac{(\eta\xi_1-\eta_1\xi)^2}{\xi\xi_1\xi_2}\right|
\]
that $|\mathscr{R}|\lesssim N^3$; so that $|\M+\ii\cchi|\lesssim N^4$. Note that for any $\zeta=(\xi,\eta)\in\rr^2$ with $\xi\in[3N/2,3N]$ and $\eta\in[-4N^2,9N^2]$, we have $|k_\zeta|\gtrsim N^3$.

Now, for $0<\epsilon\ll1$ fixed, we choose a sequence of times $\{t_N\}_N$ defined by $t_N=N^{-4-\epsilon}$. For $N\gg1$, it can be easily seen that $\ee^{-t\vrx}>C>0$. Hence
\[
\left|\ee^{t\vrx}\K\right|=\frac{1}{N^{4+\epsilon}}+O\left(\frac{1}{N^{4+2\epsilon}}\right).
\]
This implies that
\[
\left|\int_{k_\zeta}\ee^{t\vrx}\K\dd\zeta_1\right|\gtrsim|k_\zeta|N^{-4-\epsilon}\gtrsim N^{-1-\epsilon}.
\]
Therefore it follows from \eqref{norm-hs} that
\[
1\gtrsim\|u_{2,N}\|_{H^{s,0}}^2\gtrsim N^{-4s-6}N^{-2-2\epsilon}
\int_{3N/2}^{3N}
\int_{-4N^2}^{9N^2}(|\xi|+\xi^2)^2(1+\xi^2)^s\dd\zeta
\gtrsim
N^{-2s-1-2\epsilon}.
\]
This leads to a contradiction for $N\gg1$, since we have $s<-1/2-2\epsilon$; and the proof of Theorem \ref{ill-thm} is complete.
\fim
\section*{Acknowledgments}
\pdfbookmark[0]{Acknowledgments}{Acknowledgments}
The author would like to thank  Luc Molinet, Jerry L. Bona and Felipe Linares for their suggestions and comments.

\end{document}